\documentclass[amssymb, 11pt]{amsart}
\usepackage{amsmath}
\usepackage{amsthm}
\usepackage{amsfonts}
\newcommand{\pp}{\mathbb{P}}
\newtheorem{thm}{Theorem}[section]
\newtheorem{theorem}[thm]{Theorem}

\newtheorem{prop}[thm]{Proposition}

\newtheorem{lemma}[thm]{Lemma}

\usepackage{varioref}
\labelformat{section}{Section~#1}
\labelformat{figure}{Figure~#1}
\labelformat{table}{Table~#1}

\usepackage[numbers]{natbib}
\date{}

\begin{document}

\title [Mixing time for random walk on rooted trees] {A sharp analysis of the mixing time for random walk on rooted trees}

\author{Jason Fulman}

\address{Department of Mathematics, University of Southern California,
Los Angeles, CA 90089, fulman@usc.edu}



\thanks{Version of August 5, 2009.}

\begin{abstract} We define an analog of Plancherel measure for the set of
rooted unlabeled trees on $n$ vertices, and a Markov chain which has this measure
as its stationary distribution. Using the combinatorics of commutation relations,
we show that order $n^2$ steps are necessary and suffice for convergence to the
stationary distribution.
\end{abstract}

\maketitle

\section{Introduction}\label{sec1}

The Plancherel measure of the symmetric group is a probability measure on the irreducible
representations of the symmetric group which chooses a representation with probability
proportional to the square of its dimension. Equivalently, the irreducible representations
of the symmetric group are parameterized by partitions $\lambda$ of $n$, and the Plancherel
measure chooses a partition $\lambda$ with probability \begin{equation} \label{1st} \frac{n!}{\prod_{x \in
\lambda} h(x)^2} \end{equation} where the product is over boxes in the partition and $h(x)$ is the
hooklength of a box. The hooklength of a box $x$ is defined as 1 + number of boxes in
same row as x and to right of x + number of boxes in same column of x and below x. For example
we have filled in each box in the partition of 7 below with its hooklength
\[ \begin{array}{c c c c}
                \framebox{6}& \framebox{4}& \framebox{2}& \framebox{1}
                \\ \framebox{3}& \framebox{1}&& \\ \framebox{1} &&&
                \end{array}, \] and the Plancherel measure would
choose this partition with probability $\frac{7!}{(6*4*3*2)^2}$. There has been
significant interest in the statistical properties of partitions chosen from
Plancherel measure of the symmetric group; for this the reader can consult \cite{AD3},
\cite{BOO}, \cite{De} and the many references therein.

In this paper we define a similar measure on the set of rooted, unlabeled trees on $n$ vertices. We place
the root vertex on top, and the four rooted trees on 4 vertices are depicted below:

\newpage

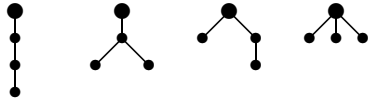
\begin{figure}[h]
\setlength{\unitlength}{.14in}
\begin{picture}(10,0)
\put(0,0){\circle*{.6}}
\put(0,0){\line(0,-1){1}}
\put(0,-1){\circle*{.4}}
\put(0,-1){\line(0,-1){1}}
\put(0,-2){\circle*{.4}}
\put(0,-2){\line(0,-1){1}}
\put(0,-3){\circle*{.4}}
\put(4,0){\circle*{.6}}
\put(4,0){\line(0,-1){1}}
\put(4,-1){\circle*{.4}}
\put(4,-1){\line(1,-1){1}}
\put(4,-1){\line(-1,-1){1}}
\put(5,-2){\circle*{.4}}
\put(3,-2){\circle*{.4}}
\put(8,0){\circle*{.6}}
\put(8,0){\line(1,-1){1}}
\put(8,0){\line(-1,-1){1}}
\put(9,-1){\circle*{.4}}
\put(7,-1){\circle*{.4}}
\put(9,-1){\line(0,-1){1}}
\put(9,-2){\circle*{.4}}
\put(12,0){\circle*{.6}}
\put(12,0){\line(-1,-1){1}}
\put(12,0){\line(0,-1){1}}
\put(12,0){\line(1,-1){1}}
\put(11,-1){\circle*{.4}}
\put(12,-1){\circle*{.4}}
\put(13,-1){\circle*{.4}}
\end{picture}
\end{figure}

\vspace{.3in}

This measure chooses a rooted tree with probability \begin{equation} \label{0}
\pi(t) = \frac{n \cdot 2^{n-1}}{|SG(t)| \prod_{v \in t} h(v)^2}, \end{equation}
where $h(v)$ is the size of the subtree with root $v$, and $|SG(t)|$ is a certain symmetry factor associated
to the tree $t$ (precise definitions are given in \ref{tree}). We do not know that this
measure has applications similar to the Plancherel measure of the symmetric group, but the resemblance is
striking. Moreover, there are Hopf algebras in the physics literature whose generators are rooted trees
(Kreimer's Hopf algebra \cite{CKr},\cite{Kr} a Hopf algebra of Connes and Moscovici \cite{CM}, and a Hopf algebra
of Grossman and Larson \cite{GL}), and as a paper of Hoffman \cite{Ho} makes clear, the combinatorics of these Hopf algebras is very close to the combinatorics we use in this paper.

In fact the main object we study is a Markov chain $K$ which has $\pi$ as its stationary distribution; this Markov
chain is defined in \ref{tree} and involves removing a single terminal vertex and reattaching it. There are
several ways of quantifying the convergence rate of a Markov chain on a state space $X$ to its stationary distribution; we use the maximal separation distance after $r$ steps, defined as
\[ s^*(r) := \max_{x,y \in X} \left[ 1 - \frac{K^r(x,y)}{\pi(y)} \right], \] where $K^r(x,y)$ is the chance
of transitioning from $x$ to $y$ after $r$ steps. In general it can be quite tricky even to determine which $x,y$
attain the maximum in the definition of $s^*(r)$. We do this, and prove that for $c>0$ fixed,
\[ \lim_{n \rightarrow \infty} s^*(cn^2) = \sum_{i=3}^{\infty} \frac{(-1)^{i-1}}{2} (2i-1)(i+1)(i-2) e^{-ci(i-1)}.\] There are very few Markov chains for which such precise asymptotics are known. Our proof method uses a commutation relation of a growth and pruning operator on rooted trees (due to Hoffman \cite{Ho}), a formula for the eigenvalues of $K$, and ideas from \cite{Fu}. Details appear in \ref{main}.

We mention that the Markov chain $K$ is very much in the spirit of the down-up chains (on the state space of partitions) studied in \cite{BO}, \cite{BO2}, \cite{Fu}, \cite{Fu3}, \cite{Pe}. There are also similarities to certain random walks on phylogenetic trees (cladograms) studied in \cite{Al}, \cite{Fo}, \cite{Sc}. Our methods only partly apply to these walks (the geometry of the two spaces of trees is different), so this will be studied in another work.

To close the introduction, we mention two reasons why it can be useful to understand a Markov chain $K$ whose stationary distribution $\pi$ is of interest. First, in analogy with Plancherel measure of the symmetric group, one can hope to use Stein's method (\cite{Fu3}) or other techniques (\cite{BO}) to study statistical properties of $\pi$. Second, convergence rates of $K$ can lead to concentration inequalities for statistics of $\pi$ \cite{C}.

\section{Background on Markov chains} \label{background}

We will be concerned with the theory of finite Markov chains. Thus $X$
will be a finite set (in our case the set of rooted unlabeled trees on $n$
vertices) and $K$ a matrix indexed by $X \times X$ whose rows
sum to 1. Let $\pi$ be a probability distribution on $X$ such that $K$
is reversible with respect to $\pi$; this means that $\pi(x) K(x,y) = \pi(y) K(y,x)$
for all $x,y$ and implies that $\pi$ is a stationary distribution for the Markov
chain corresponding to $K$ (i.e. that $\pi(x)= \sum_y \pi(y) K(y,x)$ for all
$x$).

A common way to quantify convergence rates of Markov chains is to use
separation distance, introduced by Aldous and Diaconis
\cite{AD1},\cite{AD2}. They define the separation distance of a Markov
chain $K$ started at $x$ as \[ s(r) = \max_{y} \left[ 1 -
\frac{K^r(x,y)}{\pi(y)} \right] \] and the maximal separation distance of
the Markov chain $K$ as \[ s^*(r) = \max_{x,y} \left[ 1 -
\frac{K^r(x,y)}{\pi(y)} \right]. \]

They show that the maximal separation distance has the nice properties:
\begin{itemize}
\item \[ \frac{1}{2} \max_{x} \sum_y \left| K^r(x,y)-\pi(y) \right| \leq s^*(r) \]
\item (monotonicity) $s^*(r_1) \leq s^*(r_2)$, $r_1 \geq r_2$
\item (submultiplicativity) $s^*(r_1+r_2) \leq s^*(r_1)s^*(r_2)$
\end{itemize}

\section{Combinatorics of rooted trees} \label{tree}

For a finite rooted tree $t$, we let $|t|$ denote the number of vertices of $t$;
$\mathcal{T}_n$ will be the set of rooted unlabeled trees on $n$ vertices. For example $\mathcal{T}_1
=\{ \bullet \}$ consists of only the root vertex, and the four elements of $\mathcal{T}_4$
were depicted in the introduction. Letting $T_n=|\mathcal{T}_n|$ and $T_0=0$, there is a recursion
\[ \sum_{n \geq 1} T_n \cdot x^n = x \prod_{n \geq 1} (1-x^n)^{-T_{n}} \] from which one obtains
$T_1=1, T_2=1, T_3=2, T_4=4, T_5=9, T_6=20$, etc. (see \cite{Sl} for more information on this sequence).

A rooted tree can be viewed as a directed graph by directing all edges away from the root,
and a vertex is called terminal if it has no outgoing edge. There is a
partial order $\preceq$ on the set $\mathcal{T}$ of all finite rooted trees defined by
letting $t$ be covered by $t'$ exactly when $t$ can be obtained from $t'$ by
removing a single terminal vertex and the edge into it; we denote this by $t \nearrow t'$
or $t' \searrow t$.

When $t \nearrow t'$, one can define two quantities

\[ n(t,t') = | \mbox{vertices of t to which a new edge can be added to get t}' | \] and
\[ m(t,t') = | \mbox{edges of t}' \mbox{ which when removed give t} |. \]
These need not be equal, as can be seen by taking $t,t'$ to be:

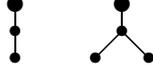
\begin{figure}[h]
\setlength{\unitlength}{.14in}
\begin{picture}(5,0)
\put(0,0){\circle*{.6}}
\put(0,0){\line(0,-1){1}}
\put(0,-1){\circle*{.4}}
\put(0,-1){\line(0,-1){1}}
\put(0,-2){\circle*{.4}}
\put(4,0){\circle*{.6}}
\put(4,0){\line(0,-1){1}}
\put(4,-1){\circle*{.4}}
\put(4,-1){\line(1,-1){1}}
\put(4,-1){\line(-1,-1){1}}
\put(5,-2){\circle*{.4}}
\put(3,-2){\circle*{.4}}
\end{picture}
\end{figure}

\vspace{.3in} Then $n(t,t')=1$ and $m(t,t')=2$.

Let $\mathbb{C} \mathcal{T}_n$ denote the complex vector space with basis the elements of $\mathcal{T}_n$.
For $n \geq 1$, Hoffman \cite{Ho} defines a growth operator ${\it G}: \mathbb{C} \mathcal{T}_n \mapsto \mathbb{C} \mathcal{T}_{n+1}$
by \[ \it{G}(t) = \sum_{t' \searrow t} n(t,t') t', \] and for $n \geq 2$ a pruning operator
${\it P}: \mathbb{C} \mathcal{T}_n \mapsto \mathbb{C} \mathcal{T}_{n-1}$ by \[ \it{P}(t)= \sum_{t' \nearrow t} m(t',t) t'.\]
One sets ${\it P} (\bullet) = 0$.

One can extend the definitions of $m(t,t')$ and $n(t,t')$ to any pair of rooted trees $t,t'$ with
$|t'|-|t|=k \geq 0$ by setting
\[ \it{G}^k(t) = \sum_{|t'|=|t|+k} n(t,t') t' \] and \[ \it{P}^k(t') = \sum_{|t|=|t'|-k} m(t,t') t.\]
Since $\bullet \preceq t$ for all $t$, one can think of $n(\bullet,t)$ as the number of ways to
build up $t$, and of $m(\bullet,t)$ as the number of ways to take $t$ apart by sequentially removing
terminal edges. To simplify notation, we let $n(t)=n(\bullet,t)$ and $m(t)=m(\bullet,t)$. For example,
the reader can check that the four trees $t_1,t_2,t_3,t_4$

\begin{figure}[h]
\setlength{\unitlength}{.14in}
\begin{picture}(10,0)
\put(0,0){\circle*{.6}}
\put(0,0){\line(0,-1){1}}
\put(0,-1){\circle*{.4}}
\put(0,-1){\line(0,-1){1}}
\put(0,-2){\circle*{.4}}
\put(0,-2){\line(0,-1){1}}
\put(0,-3){\circle*{.4}}
\put(4,0){\circle*{.6}}
\put(4,0){\line(0,-1){1}}
\put(4,-1){\circle*{.4}}
\put(4,-1){\line(1,-1){1}}
\put(4,-1){\line(-1,-1){1}}
\put(5,-2){\circle*{.4}}
\put(3,-2){\circle*{.4}}
\put(8,0){\circle*{.6}}
\put(8,0){\line(1,-1){1}}
\put(8,0){\line(-1,-1){1}}
\put(9,-1){\circle*{.4}}
\put(7,-1){\circle*{.4}}
\put(9,-1){\line(0,-1){1}}
\put(9,-2){\circle*{.4}}
\put(12,0){\circle*{.6}}
\put(12,0){\line(-1,-1){1}}
\put(12,0){\line(0,-1){1}}
\put(12,0){\line(1,-1){1}}
\put(11,-1){\circle*{.4}}
\put(12,-1){\circle*{.4}}
\put(13,-1){\circle*{.4}}
\end{picture}
\end{figure}

\vspace{.4in}

satisfy $n(t_1)=1, m(t_1)=1; n(t_2)=1, m(t_2)=2; n(t_3)=3, m(t_3)=3; n(t_4)=1, m(t_4)=6$ respectively.

There is a ``hook-length'' type formula for $m(t)$ in the literature. Namely if $t$ has $n$ vertices,
\begin{equation} \label{hook} m(t) = \frac{n!}{\prod_{v \in t} h(v)} \end{equation} where $h(v)$ is the number of vertices in the subtree with root $v$; see Section 22 of \cite{St1} or Exercise 5.1.4-20 of \cite{Kn} for a proof.

As for $n(t)$, it is also known as the ``Connes-Moscovici weight'' \cite{Kr}. To give a formula for it, we use the concept of the symmetry group $SG(t)$ of a tree. For $v$ a vertex of $T$ with children $\{v_1,\cdots,v_k\}$, $SG(t,v)$ is the group generated by the
permutations that exchange the trees with roots $v_i$ and $v_j$ when they are isomorphic rooted trees; then $SG(t)$ is defined as the direct product \[ SG(t) = \prod_{v \in T} SG(t,v).\] It is proved in \cite{Kr} that \begin{equation} \label{sym} n(t)= \frac{m(t)}{|SG(t)|}.\end{equation} More generally, Proposition 2.5 of \cite{Ho} shows that
\begin{equation} \label{sym2} n(s,t) |SG(t)| = m(s,t) |SG(s)| \end{equation} when $|s| \leq |t|$.

\vspace{.1in}

{\bf Definition 1} We define a probability measure $\pi_n$ on the set of rooted (unlabeled) trees of size $n$
by \begin{equation} \label{meas} \pi_n(t) = \frac{m(t) n(t)}{\prod_{i=2}^n {i \choose 2}} = \frac{n \cdot 2^{n-1}}{|SG(t)| \prod_{v \in t} h(v)^2}.\end{equation}

It follows from Proposition 2.8 of \cite{Ho} that $\pi$ is in fact a probability measure (i.e. that the probabilities sum to 1). The second equality in \eqref{meas} follows from equations \eqref{hook} and \eqref{sym}. The reader can check that the
four trees $t_1,t_2,t_3,t_4$

\begin{figure}[h]
\setlength{\unitlength}{.14in}
\begin{picture}(10,0)
\put(0,0){\circle*{.6}}
\put(0,0){\line(0,-1){1}}
\put(0,-1){\circle*{.4}}
\put(0,-1){\line(0,-1){1}}
\put(0,-2){\circle*{.4}}
\put(0,-2){\line(0,-1){1}}
\put(0,-3){\circle*{.4}}
\put(4,0){\circle*{.6}}
\put(4,0){\line(0,-1){1}}
\put(4,-1){\circle*{.4}}
\put(4,-1){\line(1,-1){1}}
\put(4,-1){\line(-1,-1){1}}
\put(5,-2){\circle*{.4}}
\put(3,-2){\circle*{.4}}
\put(8,0){\circle*{.6}}
\put(8,0){\line(1,-1){1}}
\put(8,0){\line(-1,-1){1}}
\put(9,-1){\circle*{.4}}
\put(7,-1){\circle*{.4}}
\put(9,-1){\line(0,-1){1}}
\put(9,-2){\circle*{.4}}
\put(12,0){\circle*{.6}}
\put(12,0){\line(-1,-1){1}}
\put(12,0){\line(0,-1){1}}
\put(12,0){\line(1,-1){1}}
\put(11,-1){\circle*{.4}}
\put(12,-1){\circle*{.4}}
\put(13,-1){\circle*{.4}}
\end{picture}
\end{figure}

\vspace{.4in}

are assigned probabilities $1/18, 1/9, 1/2, 1/3$ respectively.

\vspace{.2in}

{\bf Definition 2} We define upward transition probabilities from $t \in \mathcal{T}_{n-1}$ to $t' \in \mathcal{T}_{n}$ by \[ P_u(t,t') = \frac{m(t,t') n(t')}{{n \choose 2} n(t)} =  \frac{n(t,t') m(t')}{{n \choose 2} m(t) }\] and downward transition probabilities from $t \in \mathcal{T}_{n}$ to $t' \in \mathcal{T}_{n-1}$ by \[ P_d(t,t') = \frac{m(t',t) m(t')}{m(t)}. \]

It is clear from the definitions that the downward transition probabilities sum to 1. The second equality in the definition of $P_u$ is from \eqref{sym} and \eqref{sym2}, and it follows from Proposition 2.8 of \cite{Ho} that the upward transition probabilities sum to 1. We define a ``down-up'' Markov chain with state space $\mathcal{T}_n$ by composing the down chain with the up chain, i.e. \begin{eqnarray*} K(t,t') & = & \sum_{s \nearrow t,t'} P_d(t,s) P_u(s,t') \\
& = &  \sum_{s \nearrow t,t'} \frac{m(s,t) m(s)}{m(t)} \frac{n(s,t') m(t')}{{n \choose 2} m(s)} \\
& = & \frac{m(t')}{{n \choose 2} m(t)} \sum_{s \nearrow t,t'} m(s,t) n(s,t'). \end{eqnarray*}

Thus we deduce the crucial relation

\begin{equation} \label{cruc} K_n = \frac{1}{{n \choose 2}} A {\it G} {\it P} A^{-1}_n
\end{equation} where $A$ is the diagonal matrix which multiplies a tree $t$ by $m(t)$, and $P,G$ are the pruning
and growth operators. The subscript $n$ indicates that the chain is on trees of size $n$.

For example, ordering the four elements of $\mathcal{T}_4$ as

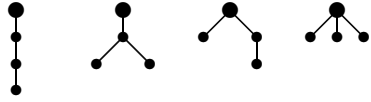
\begin{figure}[h]
\setlength{\unitlength}{.14in}
\begin{picture}(10,0)
\put(0,0){\circle*{.6}}
\put(0,0){\line(0,-1){1}}
\put(0,-1){\circle*{.4}}
\put(0,-1){\line(0,-1){1}}
\put(0,-2){\circle*{.4}}
\put(0,-2){\line(0,-1){1}}
\put(0,-3){\circle*{.4}}
\put(4,0){\circle*{.6}}
\put(4,0){\line(0,-1){1}}
\put(4,-1){\circle*{.4}}
\put(4,-1){\line(1,-1){1}}
\put(4,-1){\line(-1,-1){1}}
\put(5,-2){\circle*{.4}}
\put(3,-2){\circle*{.4}}
\put(8,0){\circle*{.6}}
\put(8,0){\line(1,-1){1}}
\put(8,0){\line(-1,-1){1}}
\put(9,-1){\circle*{.4}}
\put(7,-1){\circle*{.4}}
\put(9,-1){\line(0,-1){1}}
\put(9,-2){\circle*{.4}}
\put(12,0){\circle*{.6}}
\put(12,0){\line(-1,-1){1}}
\put(12,0){\line(0,-1){1}}
\put(12,0){\line(1,-1){1}}
\put(11,-1){\circle*{.4}}
\put(12,-1){\circle*{.4}}
\put(13,-1){\circle*{.4}}
\end{picture}
\end{figure}

\vspace{.3in}

one calculates the transition matrix \begin{equation*}
\left( K(t,t') \right) = \begin{pmatrix}
1/6 & 1/3 & 1/2 & 0 \\
1/6 & 1/3 & 1/2 & 0 \\
1/18 & 1/9 & 1/2 & 1/3 \\
0 & 0 & 1/2 & 1/2
\end{pmatrix}. \end{equation*}

The following lemma will be useful.

\begin{lemma} \label{preserve}
\begin{enumerate}
\item If $s$ is chosen from the measure $\pi_{n-1}$ and one moves from $s$ to $t$ with probability $P_u(s,t)$, then $t$ is distributed according to the measure $\pi_n$.
\item If $t$ is chosen from the measure $\pi_{n+1}$ and one moves from $t$ to $s$ with probability $P_d(t,s)$, then $s$ is distributed according to the measure $\pi_n$.
\item The ``down-up'' Markov chain $K_n$ on rooted trees of size $n$ is reversible with respect to $\pi_n$.
\end{enumerate}
\end{lemma}

\begin{proof} For part 1, one calculates that
\begin{eqnarray*} \sum_{s \nearrow t} \pi_{n-1}(s) P_u(s,t) & = & \sum_{s \nearrow t} \frac{m(s) n(s)}
{\prod_{i=2}^{n-1} {i \choose 2}} \frac{n(s,t) m(t)}{m(s) {n \choose 2}} \\
& = & \frac{m(t)}{\prod_{i=2}^n {i \choose 2}} \sum_{s \nearrow t} n(s) n(s,t) \\
& = & \frac{m(t)n(t)}{\prod_{i=2}^n {i \choose 2}} = \pi_n(t). \end{eqnarray*}

For part 2, one computes that \begin{eqnarray*}
\sum_{t \searrow s} \pi_{n+1}(t)P_d(t,s) & = & \sum_{t \searrow s} \frac{m(t) n(t)}{\prod_{i=2}^{n+1} {i \choose 2}}
\frac{m(s,t) m(s)}{m(t)} \\
& = & \frac{m(s)}{\prod_{i=2}^{n+1} {i \choose 2}} \sum_{t \searrow s} n(t) m(s,t) \\
& = & \frac{m(s)n(s)} {\prod_{i=2}^{n} {i \choose 2}} = \pi_n(s), \end{eqnarray*} where the last line follows since
the upward transition probabilities from $s$ sum to $1$.

For part 3, one calculates that \begin{eqnarray*}
\pi_n(t)K(t,t') & = & \frac{n(t) m(t')}{{n \choose 2} \prod_{i=2}^n {i \choose 2}} \sum_{s \nearrow t,t'}
m(s,t) n(s,t') \\
& = & \frac{n(t) m(t')}{{n \choose 2} \prod_{i=2}^n {i \choose 2}} |SG(t)| \sum_{s \nearrow t,t'}
\frac{n(s,t) n(s,t')}{|SG(s)|} \\
& = & \frac{n(t) m(t')}{{n \choose 2} \prod_{i=2}^n {i \choose 2}} |SG(t)| \sum_{s \nearrow t,t'}
\frac{n(s,t) m(s,t')}{|SG(t')|} \\
& = & \frac{n(t')m(t)} {{n \choose 2} \prod_{i=2}^n {i \choose 2}} \sum_{s \nearrow t,t'}
n(s,t) m(s,t') \\
& = & \pi_n(t') K(t',t). \end{eqnarray*} Note that equation \eqref{sym2} was used in equalities 2 and 3 and that equation \eqref{sym}
was used in the fourth equality.
\end{proof}

The final combinatorial fact we will need about rooted trees is the following commutation relation between the growth and pruning operators (Proposition 2.2 of \cite{Ho}) :
\begin{equation} \label{1} {\it P} {\it G}_n - {\it G} {\it P}_n = n I,
\end{equation} for all $n \geq 1$. Here $I$ is the identity operator, so the right hand side multiplies a tree by its size.

\section{Proof of main results} \label{main}

The purpose of this section is to obtain precise asymptotics for the maximal separation distance $s^*(r)$ of the Markov chain $K$ after $r$ iterations. To do this we use equation \eqref{cruc}, the commutation relation \eqref{1}, and the methodology of \cite{Fu}. To begin we determine the eigenvalues of the Markov chain $K$. The multiplicities involve the numbers $T_i$ of rooted unlabeled trees of size $i$, discussed in \ref{tree}.

\begin{prop} \label{eigenval} The eigenvalues of the Markov chain $K$ are:
$$ \begin{array}{ll} 1 & \mbox{multiplicity} \ 1 \\
1 - \frac{{i \choose 2}}{{n \choose 2}} & \mbox{multiplicity} \
T_{i}-T_{i-1} \ (3 \leq i \leq n)
\end{array} 
$$
\end{prop}

\begin{proof} Since $K_n = \frac{1}{{n \choose 2}} A {\it G} {\it P} A^{-1}_n$, it suffices to determine the eigenvalues of ${\it G} {\it P}_n$; these follow from the commutation relation \eqref{1} and Theorem 2.6 of
\cite{St2}. \end{proof}

Recall that our interest is in studying the behavior of \[ s^*(r) = \max_{t,t'} \left[ 1 - \frac{K^r(t,t')}{\pi(t')} \right].\] Proposition \ref{obtained} determines the pairs $(t,t')$ where this maximum is obtained.

\begin{prop} \label{obtained} For all values of $r$, the quantity $1 - \frac{K^r(t,t')}{\pi(t')}$ is maximized by letting $t$ be the unique rooted tree with one terminal vertex and $t'$ be the unique tree with $n-1$ terminal vertices, or by letting $t'$ be the unique rooted tree with one terminal vertex and $t$ be the unique tree with $n-1$ terminal vertices.
\end{prop}

For instance when $n=5$ the two relevant trees are

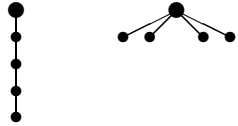
\begin{figure}[h]
\setlength{\unitlength}{.14in}
\begin{picture}(10,0)
\put(0,0){\circle*{.6}}
\put(0,0){\line(0,-1){1}}
\put(0,-1){\circle*{.4}}
\put(0,-1){\line(0,-1){1}}
\put(0,-2){\circle*{.4}}
\put(0,-2){\line(0,-1){1}}
\put(0,-3){\circle*{.4}}
\put(0,-3){\line(0,-1){1}}
\put(0,-4){\circle*{.4}}
\put(6,0){\circle*{.6}}
\put(6,0){\line(-1,-1){1}}
\put(6,0){\line(-2,-1){2}}
\put(6,0){\line(1,-1){1}}
\put(6,0){\line(2,-1){2}}
\put(5,-1){\circle*{.4}}
\put(4,-1){\circle*{.4}}
\put(7,-1){\circle*{.4}}
\put(8,-1){\circle*{.4}}
\end{picture}
\end{figure}

\vspace{.3in}

\begin{proof} By relation \eqref{cruc}, we seek the $t,t'$ minimizing \[ \frac{K^r(t,t')}{\pi(t')} = \frac{m(t') ({\it G} {\it P})^r_{n} [t,t']}{{n \choose 2}^r m(t) \pi(t')}.\] By the commutation relation \eqref{1} and Proposition 4.5 of \cite{Fu}, \[  ({\it G} {\it P})^r_{n} = \sum_{k=0}^{n} A_{n}(r,k) {\it G}^k {\it P}^k_{n} \] where the $A_{n}(r,k)$ solve the recurrence \[ A_{n}(r,k) = A_{n}(r-1,k-1) + A_{n}(r-1,k) \left[ {n \choose 2} - {n-k \choose 2} \right] \] with initial conditions $A_{n}(0,0)=1$, $A_{n}(0,m)=0$ for $m \neq 0$. Thus \begin{equation} \label{2} \frac{K^r(t,t')}{\pi(t')} = \frac{m(t')\sum_{k=0}^{n} A_{n}(r,k) {\it G}^k {\it P}^k_{n}[t,t']}{{n \choose 2}^r m(t) \pi(t')}.\end{equation}

The proposition now follows from three observations:

\begin{itemize}
\item All terms in \eqref{2} are non-negative. Indeed, this is clear from the recurrence for $A_{n}(r,k)$.

\item If $t$ is the unique rooted tree with one terminal vertex and $t'$ is the unique rooted tree with $n-1$ terminal vertices (or the same holds with $t,t'$ swapped), then the summands in \eqref{2} for $0 \leq k \leq n-3$ all vanish. Indeed, in order to move from $t$ to $t'$ by pruning $k$ vertices and then reattaching them, one must prune at least $n-2$ vertices.

\item The $k=n-2$ and $k=n-1$ summands in \eqref{2} are independent of $t,t'$. Indeed, for the $k=n-1$ summand, one has that
\begin{eqnarray*}
& &  \frac{m(t') A_{n}(r,n-1) {\it G}^{n-1} {\it P}^{n-1}_{n}[t,t']}{{n \choose 2}^r m(t) \pi(t')} \\
& = & \frac{m(t') A_{n}(r,n-1) {\it G}^{n-1} [\bullet,t']}{{n \choose 2}^r \pi(t')} \\
& = & \frac{m(t') A_{n}(r,n-1) n(t')}{{n \choose 2}^r \pi(t')} \\
& = & A_{n}(r,n-1) \frac{\prod_{i=2}^n {i \choose 2}}{{n \choose 2}^r}.
\end{eqnarray*} A similar argument shows that the $k=n-2$ summand is equal to $A_{n}(r,n-2) \frac{\prod_{i=2}^n {i \choose 2}}{{n \choose 2}^r}$.

\end{itemize}
\end{proof}

{\it Remark:} The proof of Proposition \ref{obtained} shows that
\[ s^*(r) = 1- \frac{\prod_{i=2}^n {i \choose 2}}{{n \choose 2}^r} \left[ A_{n}(r,n-2) + A_{n}(r,n-1) \right],\] where $A_{n}(r,k)$ is the solution to the recurrence in the proof of Proposition \ref{obtained}.

\vspace{.2in}

In Theorem \ref{asymp}, we give an explicit formula for $s^*(r)$ and determine its asymptotic behavior.

\begin{theorem} \label{asymp} Let $s^*(r)$ be the maximal separation distance after $r$ iterations of the down-up Markov chain $K$ on the space of rooted trees on $n$ vertices.
\begin{enumerate}
\item For $r \geq 1$, $s^*(r)$ is equal to
\[ \sum_{i=3}^{n-1} (-1)^{i-1} \frac{(2i-1)(i+1)(i-2)(n!)^2}{2n (n-i)! (n+i-1)!} \left( 1 - \frac{{i \choose 2}}{{n \choose 2}} \right)^r .\]
\item For $c>0$ fixed,
\[ \lim_{n \rightarrow \infty} s^*(cn^2) = \sum_{i=3}^{\infty} \frac{(-1)^{i-1}}{2} (2i-1)(i+1)(i-2) e^{-ci(i-1)}.\]
\end{enumerate}
\end{theorem}

\begin{proof} By Proposition \ref{obtained}, the maximal separation distance is attained when $t$ is the unique rooted tree with one terminal vertex and $t'$ is the unique rooted tree with $n-1$ terminal vertices. Note that it takes $n-2$ iterations of the Markov chain $K$ to move from $t$ to $t'$. By Proposition \ref{eigenval}, $K$ has $n-1$ distinct eigenvalues (one more than the Markov chain distance between $t$ and $t'$), so it follows from Proposition 5.1 of \cite{Fu2} that \begin{equation} \label{sep} s^*(r) = \sum_{i=3}^{n} \lambda_i^r \left[ \prod_{j \neq i} \frac{1-\lambda_j}{\lambda_i-\lambda_j} \right], \end{equation} where $1$, $\lambda_i=1-\frac{{i \choose 2}}{{n \choose 2}}$, $i=3,\cdots,n$ are the distinct eigenvalues of $K$. For $r \geq 1$, this is equal to \begin{equation} \begin{split} \label{11}  &  \sum_{i=3}^{n-1} \left( 1-\frac{{i \choose 2}}{{n \choose 2}} \right)^r \prod_{j \neq i \atop 3 \leq j \leq n}
\frac{{j \choose 2}}{{j \choose 2} - {i \choose 2}} \\
& =  \sum_{i=3}^{n-1} \left( 1-\frac{{i \choose 2}}{{n \choose 2}} \right)^r \prod_{j \neq i \atop 3 \leq j \leq n}
\frac{j(j-1)}{(j+i-1)(j-i)}, \end{split} \end{equation} and the first assertion follows by elementary simplifications.

For part 2 of the theorem, it is enough to show that for $c>0$ fixed, there is a constant $i_c$ (depending on $c$ but not on $n$) such that for $i \geq i_c$, the summands in part 1 of the theorem are decreasing in magnitude (and alternating in sign). Part 2 follows from this claim, since then one can take limits for each fixed $i$. For $i \geq 2 \sqrt{n}$ one checks that \[ \frac{(2i-1)(i+1)(i-2)(n!)^2}{2n (n-i)! (n+i-1)!} \] is a decreasing function of $i$. To handle the case of $i \leq 2 \sqrt{n}$, one need only show that \begin{equation} \label{need} \frac{(n-i)(2i+1)(i+2)(i-1)}{(n+i)(2i-1)(i+1)(i-2)} \frac{\exp(cn^2 \log(1-{i+1 \choose 2}/{n \choose 2}))} {\exp(cn^2 \log(1-{i \choose 2}/{n \choose 2}))} < 1 \end{equation} for $i \geq i_c$, a constant depending on $c$ but not on $n$. This is easily established, since using the inequalities $\log(1-x) \leq -x$ for $x>0$ in the numerator and $\log(1-x) \geq -x-x^2$ for $0<x<\frac{1}{2}$ in the denominator gives that \[ \frac{\exp(cn^2 \log(1-{i+1 \choose 2}/{n \choose 2}))} {\exp(cn^2 \log(1-{i \choose 2}/{n \choose 2}))} \leq \exp \left[ \frac{-cn^2}{{n \choose 2}} \left( i - \frac{{i \choose 2}^2}{{n \choose 2}} \right) \right], \] and \eqref{need} follows as $i \leq 2 \sqrt{n}$.
\end{proof}

Some authors who work on Markov chains similar to that studied here but on different state spaces (e.g. \cite{BO2}, \cite{Pe}) prefer to work with up-down chains instead of down-up chains. Proposition \ref{upd} shows the study of maximal separation for these two chains to be equivalent.

\begin{prop} \label{upd} Let $s^*_{UD_n}(r)$ denote the maximal separation distance after $r$ iterations of the down-up chain on $\mathcal{T}_n$, and let $s^*_{DU_n}(r)$ be the corresponding quantity for the up-down chain. Then \[ s^*_{DU_n}(r) =
s^*_{UD_{n+1}}(r+1) \] for all $n,r \geq 1$. \end{prop}

\begin{proof} An argument similar to that used to prove equation \eqref{cruc} gives that \begin{equation} \label{rel2} DU_n = \frac{1}{{n+1 \choose 2}} APGA_n^{-1} \end{equation} where $A$ is the diagonal matrix which multiplies a tree by $m(t)$, and $P,G$ are the pruning and growth operators. Combining this with the commutation relation \eqref{1}, it follows that \begin{eqnarray*} (DU)_n^r & = & \frac{1}{{n+1 \choose 2}^r} [A (nI+GP)A_n^{-1}]^r \\
& = & \frac{1}{{n+1 \choose 2}^r} \sum_{l=0}^r {r \choose l} n^{r-l} A (GP)^l A_n^{-1}. \end{eqnarray*} Arguing as in Proposition \ref{obtained}, one concludes that the same $t,t'$ maximize the separation distance. Moreover, one sees from \eqref{rel2}, commutation relation \eqref{1}, and Proposition \ref{eigenval} that the distinct eigenvalues of the up-down chain on trees of size $n$ are $1$ and $\mu_i=1 - \frac{{i \choose 2}}{{n+1 \choose 2}}$, $i=3, \cdots, n$. Thus the argument of Theorem \ref{asymp} gives that
\[ s^*_{DU_n}(r) = \sum_{i=3}^n \left( 1 - \frac{{i \choose 2}}{{n+1 \choose 2}} \right)^r \prod_{j \neq i \atop 3 \leq j \leq n} \frac{{j \choose 2}}{{j \choose 2}-{i \choose 2}}.\] The proposition now follows by making the replacements $r \rightarrow r+1$ and $n \rightarrow n+1$ in the left hand side of equation \eqref{11}. \end{proof}

To close, we note the following probabilistic interpretation of $s^*(r)$. We use the convention that a random variable $X$ is called geometric with parameter (probability of success) $p$ if $\pp(X=n)=p(1-p)^{n-1}$ for all $n \geq 1$.

\begin{prop} \label{prob} Letting $s^*(r)$ be as in Theorem \ref{asymp}, one has that $s^*(r)=\pp(T>r)$, where $T = \sum_{i=3}^{n} X_i$, and the $X_i$'s are independent geometrics with parameters $\frac{{i \choose 2}}{{n \choose 2}}$.
\end{prop}

\begin{proof} This is immediate from equation \eqref{sep} and Proposition 2.4 of \cite{Fu}. \end{proof}

We remark that representations of separation distance similar to that in Proposition \ref{prob} are in the literature for stochastically monotone birth-death chains with non-negative eigenvalues (\cite {DF}, \cite{DS}) and for some random walks on partitions \cite{Fu}. Of course the Markov chain $K$ studied in this paper is not one-dimensional.

\section*{Acknowledgments}
The author was supported by NSA grant H98230-08-1-0133 and NSF grant DMS 0802082. We thank Persi Diaconis for pointers to the literature.

\end{document}